\documentclass[12pt, reqno]{amsart}

\author[M.~Balcerzak]{Marek Balcerzak}
\address{Institute of Mathematics, Lodz University of Technology, ul. W\'{o}lcza\'{n}ska 215, 93-005 Lodz, Poland} 
\email{marek.balcerzak@p.lodz.pl}

\author[P.~Leonetti]{Paolo Leonetti}
\address{Department of Statistics, Universit\`a Bocconi, via Roentgen 1, Milan 20136, Italy}
\email{leonetti.paolo@gmail.com}

\keywords{Summable ideal; Fubini sum; convergent series.}
\subjclass[2010]{Primary: 40A05. Secondary: 54A20.}

\title{Convergent subseries of divergent series}

\usepackage[T1]{fontenc}
\usepackage{amsmath}
\usepackage{amssymb}
\usepackage{amsthm}
\usepackage[left=2.5cm, right=2.5cm, top=3cm]{geometry}
\usepackage{hyperref}
\usepackage{fancyhdr}
\usepackage{enumitem}
\usepackage{comment}
\usepackage{nicefrac}
\usepackage{bm}
\usepackage{mathrsfs}
\usepackage{graphicx}
\usepackage[utf8]{inputenc}
\usepackage{cancel}
\usepackage{mathtools}

\newcommand{\vertiii}[1]{{\left\vert\kern-0.25ex\left\vert\kern-0.25ex\left\vert #1 
    \right\vert\kern-0.25ex\right\vert\kern-0.25ex\right\vert}}
		
\AtBeginDocument{%
   \def\MR#1{}
}

\newtheorem{thm}{Theorem}[section]

\newtheorem{prop}[thm]{Proposition}

\theoremstyle{definition} 
\newtheorem{defi}[thm]{Definition}
\let\olddefi\defi
\renewcommand{\defi}{\olddefi\normalfont}
\newtheorem{example}[thm]{Example}
\let\oldexample\example
\renewcommand{\example}{\oldexample\normalfont}
\newtheorem{question}[thm]{Question}
\newtheorem{rmk}[thm]{Remark}
\let\oldrmk\rmk
\renewcommand{\rmk}{\oldrmk\normalfont}

\hypersetup{
    pdftitle={},
    pdfauthor={},
    pdfmenubar=false,
    pdffitwindow=true,
    pdfstartview=FitH,
    colorlinks=true,
    linkcolor=blue,
    citecolor=green,
    urlcolor=cyan
}

\providecommand{\MR}[1]{}

\providecommand{\MR}{\relax\ifhmode\unskip\space\fi MR }

\providecommand{\href}[2]{#2}

\begin{document}

\begin{abstract}
Let $\mathscr{X}$ be the set of positive real sequences $x=(x_n)$ such that the series $\sum_n x_n$ is divergent. For each $x \in \mathscr{X}$, let $\mathcal{I}_x$ be the collection of all $A\subseteq \mathbf{N}$ such that the subseries $\sum_{n \in A}x_n$ is convergent. 
Moreover, let $\mathscr{A}$ be the set of sequences $x \in \mathscr{X}$ such that $\lim_n x_n=0$ and $\mathcal{I}_x\neq \mathcal{I}_y$ for all sequences $y=(y_n) \in \mathscr{X}$ with $\liminf_n y_{n+1}/y_n>0$. 
We show that $\mathscr{A}$ is comeager and that contains uncountably many sequences $x$ which generate pairwise nonisomorphic ideals $\mathcal{I}_x$. 
This answers, in particular, an open question recently posed by M. Filipczak and G. Horbaczewska. 
\end{abstract}

\maketitle
\thispagestyle{empty}


\section{Introduction}

Let $\mathscr{X}$ be the set of positive real sequences $x=(x_n)$ with divergent series $\sum_n x_n$. 
For each $x \in \mathscr{X}$, let $\mathcal{I}_x$ the collection of sets of positive integers $A$ such that the (possibly finite) subseries indexed by $A$ is convergent, that is, 
\begin{equation}\label{eq:summabledefinition}
\mathcal{I}_x:=\left\{A\subseteq \mathbf{N}: \sum_{n \in A}x_n < \infty\right\}.
\end{equation}

Note that each $\mathcal{I}_x$ is closed under finite unions and subsets, i.e., it is an \emph{ideal}. 
Moreover, it contains the collection $\mathrm{Fin}$ of finite sets $A\subseteq \mathbf{N}$, and and it is different from the power set $\mathcal{P}(\mathbf{N})$. 
Following \cite{MR1711328}, a collection of sets of the type \eqref{eq:summabledefinition} is called \emph{summable ideal}. 
It is not difficult to see that every infinite set of $\mathbf{N}$ contains an infinite subset in $\mathcal{I}_x$ if and only if $\lim_n x_n=0$. 
Accordingly, define 
$$
\mathscr{Z}:=\mathscr{X} \cap c_0=\{x \in \mathscr{X}: \lim_{n\to \infty} x_n=0\}.
$$

It is known that the families $\mathcal{I}_x$ defined in \eqref{eq:summabledefinition} are ``small'', both in the measure-theoretic sense and the categorical sense, meaning that ``almost all'' subseries diverge, see \cite{MR0009997, MR2189797, MR0493027, MR0179507}. 
Related results in the context of filter convergence have been given in \cite{MR3712964, BPW2018, 
LMM}. 
The set of limits of convergent subseries of a given series $\sum_n x_n$, which is usually called  ``achievement set'', has been studied in \cite{
MR2812282, MR3910632, MR3418208}. 
Of special interest have been specific subseries 
of the harmonic series $\sum_n \frac{1}{n}$; see, e.g., 
\cite{
MR3779223, 
MR2416253, 
MR3372065, 
MR508230}.

Roughly, the question that we are going to answer is the following: Is it true that for each $x \in \mathscr{Z}$ there exists $y \in \mathscr{X}$ such that $\mathcal{I}_x=\mathcal{I}_y$ and $y_n$ ``does not oscillates too much''? 

Hoping for a characterization of the class of summable ideals $\mathcal{I}_x$ with $x \in \mathscr{Z}$, M. Filipczak and G. Horbaczewska asked recently in \cite{Filipczak} the following:
\begin{question}\label{question:filipczak}
Is it true that for each $x \in \mathscr{Z}$ there exists $y \in \mathscr{X}$ such that $\mathcal{I}_x=\mathcal{I}_y$ and 
$$
\forall n\in \mathbf{N},\quad 
\frac{y_{n+1}}{y_n} \ge \frac{n}{n+2}\,\, ?
$$
\end{question}

We show in Theorem \ref{thm:existencefubini} below that the answer is negative in a strong sense. 
To this aim, define 
$$
\mathscr{Y}:=\left\{y \in \mathscr{X}: \liminf_{n\to \infty}\, \frac{y_{n+1}}{y_n} >0\right\},
$$
and let $\sim$ be the equivalence relation on $\mathscr{X}$ so that two sequences are identified if they generate the same ideal, so that
$$
\forall x,y \in \mathscr{X}, \quad 
x\sim y 
\Longleftrightarrow \left(\forall A\subseteq \mathbf{N}, \,\,\sum_{n \in A}x_n<\infty \Longleftrightarrow \sum_{n \in A}y_n<\infty\right).
$$

First, we show that the set of pairs $(x,y) \in \mathscr{X}^2$ such that $x$ is $\sim$-equivalent to $y$ is topologically well behaved:
\begin{prop}\label{prop:coanalytic}
$\sim$ is a coanalytic relation on $\mathscr{X}$.
\end{prop}

Then, we answer Question \ref{question:filipczak} by showing that:
\begin{thm}\label{thm:existencefubini}
There exists $x \in \mathscr{Z}$ such that $x\not\sim y$ for all $y \in \mathscr{Y}$.
\end{thm}

In light of the explicit example which will be given in the proof of Theorem \ref{thm:existencefubini}, one may ask about the topological largeness of the set of such sequences. To be precise, is it true that 
$$
\mathscr{A}:=\left\{x \in \mathscr{Z}: \forall y \in \mathscr{Y}, x\not\sim y\right\}
$$
is a set of second Baire category, i.e., not topologically small? Note that the question is really meaningful since $\mathscr{Z}$ is completely metrizable (hence by Baire's category theorem $\mathscr{Z}$ is not meager in itself): this follows by Alexandrov's theorem 
\cite[Theorem 3.11]{MR1321597} 
and the fact that
$$
\mathscr{Z}=\bigcap_{n\ge 1}\left\{x\in c_0: x_n>0\right\} \cap \bigcap_{m\ge 1}\bigcup_{k\ge 1}\left\{x \in c_0: x_1+\cdots+x_k>m\right\}
$$
is a $G_\delta$-subset of the Polish space $c_0$. With the premises, we show that $\mathscr{A}$ is comeager, that is, $\mathscr{Z}\setminus \mathscr{A}$ is a set of first Baire category:
\begin{thm}\label{thm:marekmeagercounterexamples}
$\mathscr{A}$ is comeager in $\mathscr{Z}$. In particular, $\mathscr{A}$ is uncountable. 
\end{thm}

We remark that Theorem \ref{thm:marekmeagercounterexamples} gives an additional information on relation $\sim$. 
Since it is a coanalytic equivalence relation by Proposition \ref{prop:coanalytic}, we can appeal to the deep theorem by Silver \cite[Theorem 35.20]{MR1321597} which states that every coanalytic equivalence relation on a Polish space either has countably many equivalence classes or there is a perfect set consisting of non-equivalent pairs. 
Thanks to Theorem \ref{thm:marekmeagercounterexamples}, the latter holds for the relation $\sim$ in a strong form. 
Indeed, every pair in $\mathscr{A}\times \mathscr{Y}$ does not belong to $\sim$, where $\mathscr{A}$ is comeager (hence it contains a $G_\delta$-comeager subset) and $\mathscr{Y}$ is an uncountable $F_\sigma$-set. 
Therefore $\mathscr{A}\times \mathscr{Y}$ contains a product of two perfect sets by \cite[Theorem 13.6]{MR1321597}.


Lastly, on a similar direction, we strenghten the fact that $\mathscr{A}$ is uncountable by proving that exist uncountably many sequences in $\mathscr{A}$ which generate pairwise nonisomorphic ideals (here, recall that two ideals $\mathcal{I}, \mathcal{J}$ are isomorphic if there exists a bijection $f: \mathbf{N}\to \mathbf{N}$ such that $f[A] \in \mathcal{I}$ if and only if $A \in \mathcal{J}$ for all $A\subseteq \mathbf{N}$).
\begin{thm}\label{thm:examples}
There are $\mathfrak{c}$ 
sequences in $\mathscr{A}$ which generate pairwise nonisomorphic ideals.
\end{thm}

Hereafter, we use the convention that $\sum_{n\ge 1}a_n \ll \sum_{n\ge 1}b_n$, with each $a_n,b_n>0$, is a shorthand for the existence of 
$C>0$ such that $\sum_{n\le k}a_n \le C\sum_{n\le k}b_n$ for all $k \in \mathbf{N}$.

\section{Proof of Proposition \ref{prop:coanalytic}}
Equiavalently, we have to show that the set 
$
E:=\left\{(x,y) \in \mathscr{X}^2: x\not\sim y\right\}
$ 
is analytic in $\mathscr{X}^2$. For, note that $E$ is the projection on $\mathscr{X}^2$ of $E_1 \cup E_2$, where 
$$
E_1:=\left\{(A,x,y) \in \mathcal{P}(\mathbf{N}) \times \mathscr{X}^2: A \in \mathcal{I}_x\setminus \mathcal{I}_y\right\}
$$
and, similarly, 
$$
E_2:=\left\{(A,x,y) \in \mathcal{P}(\mathbf{N}) \times \mathscr{X}^2: A \in \mathcal{I}_y\setminus \mathcal{I}_x\right\}.
$$
Now, for each $n \in \mathbf{N}$, define the functions $\alpha_n, \beta_n: \mathcal{P}(\mathbf{N}) \times \mathscr{X}^2\to \mathbf{R}$ by 
$\alpha_n(A,x,y)=\sum x_k$ and $\beta_n(A,x,y)=\sum y_k$, where each sum is extended over all $k \in A$ such that $k\le n$. 
Since they are continuous, the set $(\alpha_n \le k):=\{(A,x,y)\in \mathcal{P}(\mathbf{N}) \times \mathscr{X}^2: \alpha_n(A,x,y) \le k\}$ is closed and $(\beta_n>k)$ is open for all $n,k \in \mathbf{N}$. Therefore
$$
E_1=\left(\bigcup_{k\ge 1}\bigcap_{n\ge 1}(\alpha_n \le k)\right) \cap \left(\bigcap_{k\ge 1}\bigcup_{n\ge 1}(\beta_n > k)\right)
$$
is the intersection of an $F_\sigma$-set and a $G_\delta$-set, hence it is Borel. Analogously, $E_2$ is Borel. This proves that $E$ is analytic subset of $\mathscr{X}^2$.

\section{Proof of Theorem \ref{thm:existencefubini}}
Define the sequence $x=(x_n)$ so that $x_n=\frac{1}{n}$ if $n$ is even and $x_n=\frac{1}{n\log(n+1)}$ if $n$ is odd. Note that $\lim_{n}x_n=0$ and that $\sum_n x_n=\infty$, hence $x \in \mathscr{Z}$. 
At this point, fix $y \in \mathscr{Y}$ such that $\kappa:=\liminf_n y_{n+1}/y_n>0$ and let us show that $\mathcal{I}_x \neq \mathcal{I}_y$. 

Let $\mathbf{P}$ be the set of prime numbers, with increasing enumeration $(p_n)$. By the prime number theorem we have $p_n$ is asymptotically equal to $n\log(n)$ as $n\to \infty$, hence
$$
\sum_{n \in \mathbf{P}}x_n=\sum_{n\ge 1} x_{p_n} \ll \sum_{n\ge 1} \frac{1}{p_n\log(p_n)}\ll \sum_{n\ge 2} \frac{1}{n\log^2(n)}<\infty,
$$
with the consequence that $\mathbf{P} \in \mathcal{I}_x$. In addition, $\mathbf{P}-1 \notin \mathcal{I}_x$ because
$$
\sum_{n \in \mathbf{P}-1}x_n=\sum_{n\ge 1} x_{p_n-1}\gg \sum_{n\ge 1} \frac{1}{p_n} \gg \sum_{n\ge 2}\frac{1}{n\log(n)}=\infty.
$$

Lastly, suppose for the sake of contradiction that $\mathcal{I}_x=\mathcal{I}_y$. Then we should have that $\mathbf{P} \in \mathcal{I}_y$ and, at the same time, $\mathbf{P}-1\notin \mathcal{I}_y$. The latter means that
$$
\sum_{n\ge 1}y_{p_n-1}=\infty.
$$
However, this implies that
\begin{equation}\label{eq:impossible}
\sum_{n\ge 1}y_{p_n} \gg \sum_{n\ge 1} \kappa \, y_{p_n-1}=\infty,
\end{equation}
contradicting that $\mathbf{P} \in \mathcal{I}_y$.

\section{Proof of Theorem \ref{thm:marekmeagercounterexamples}}

Consider the Banach--Mazur game defined as follows: 
Players I and II choose alternatively nonempty open subsets of $\mathscr{Z}$ as a nonincreasing chain $U_1 \supseteq V_1 \supseteq U_2 \supseteq V_2 \supseteq \cdots$, where Player I chooses the sets $U_m$. Player II has a winning strategy if $\bigcap_m V_m\subseteq \mathscr{A}$. Thanks to \cite[Theorem 8.33]{MR1321597}, Player II has a winning strategy if and only if $\mathscr{A}$ is comeager. Hence, the rest of the proof consists in showing that Player II has a winning strategy.

Note that the open neighborhood of a sequence $x \in \mathscr{Z}$ with radius $\varepsilon>0$ satisfies
$$
B_\varepsilon(x):=\{y \in \mathscr{Z}: \|x-y\|<\varepsilon\} \supseteq 
\left\{y \in \mathscr{Z}: \forall n \in \mathbf{N}, \, |x_n-y_n| < \nicefrac{\varepsilon}{2}\right\}.
$$
Since $x \in \mathscr{Z}\subseteq c_0$, there exists $k_0=k_0(x,\varepsilon) \in \mathbf{N}$ such that $x_n<\nicefrac{\varepsilon}{2}$ for all $n\ge k_0$. Hence
\begin{equation}\label{eq:openneighborhood}
B_\varepsilon(x) \supseteq W_\varepsilon(x):=\left\{y \in \mathscr{Z}: \forall n \ge k_0(x,\varepsilon), \, y_n < \nicefrac{\varepsilon}{2} \text{ and }\forall n<k_0(x,\varepsilon),|x_n-y_n|<\nicefrac{\varepsilon}{2}\right\}.
\end{equation}
For each $m \in \mathbf{N}$, suppose that the nonempty open set $U_m$ has been fixed by Player I. Hence, $U_m$ contains an open ball $B_{\varepsilon_m}(x^{(m)})$, for some $x^{(m)} \in \mathscr{Z}$ and $\varepsilon_m>0$. In particular, thanks to \eqref{eq:openneighborhood}, there exists a sufficiently large integer $k_0=k_0(x^{(m)}, \varepsilon_m) \in \mathbf{N}$ such that $y_n<\nicefrac{\varepsilon}{2}$ for all $y \in W_{\varepsilon_m}(x^{(m)})$ and $n \ge k_0$. Without loss of generality, let us suppose that $k_0$ is even.

At this point, let $x^\star$ be the sequence in $\mathscr{A}$ defined in the proof of Theorem \ref{thm:existencefubini}. Then, for each $m \in \mathbf{N}$, let $t_m$ be an integer such that $\max\{x^\star_{p_m},x^\star_{p_m-1}\}<t_m\cdot \frac{\varepsilon_m}{2}$ (we recall that $p_m$ stands for the $m$th prime number), and define the positive real
$$
\delta_m:=\min\left\{\frac{1}{m^2t_m}, \frac{\varepsilon_m}{2}-\frac{\max\{x^\star_{p_m},x^\star_{p_m-1}\}}{t_m}\right\}.
$$
Note that $\lim_m \delta_m=0$. Now, define the set $I_m=\mathbf{N}\cap [k_0(x^{(m)}, \varepsilon_m), k_0(x^{(m)}, \varepsilon_m)+2t_m)$ and let $z^{(m)}$ be the sequence such that 
\begin{displaymath}
\forall n \in \mathbf{N} \quad 
z^{(m)}_n=
\begin{cases}
x^\star_{p_m-1}/t_m \,\,\,\,\,\, & \text{ if }n \in I_m \text{ and } n \text{ even},\\
x^\star_{p_m}/t_m  \,\,\,\,\,\, & \text{ if }n \in I_m \text{ and } n \text{ odd},\\
x^{(m)}_n  \,\,\,\,\,\, & \text{ if }n \notin I_m.\\
\end{cases}
\end{displaymath}
Lastly, set $V_m:=B_{\delta_m}(z^{(m)})$ and note that by construction
$$
\forall m \in \mathbf{N}, \quad 
V_m \subseteq W_{\varepsilon_m}(x^{(m)}) \subseteq B_{\varepsilon_m}(x^{(m)}) \subseteq U_m,
$$
hence $V_m$ is a nonempty open set contained in $U_m$. In addition, the sequence of centers $(z^{(m)})$ is a Cauchy sequence in the complete metric space $\mathscr{Z}$. Hence it is convergent to some $z \in \mathscr{Z}$ and it is straighforward to see that $\{z\}=\bigcap_m V_m$.

To complete the proof, we need to show that $z \in \mathscr{A}$. Set $A:=\left(\bigcup_m I_m\right) \setminus 2\mathbf{N}$. Proceeding as in the proof of Theorem \ref{thm:existencefubini}, we see that 
\begin{displaymath}
\sum_{n \in A}z_n=\sum_{m\ge 1}\sum_{n \in I_m \setminus 2\mathbf{N}}z_n \le \sum_{m\ge 1}\sum_{n \in I_m \setminus 2\mathbf{N}}(z_n^{(m)}+\delta_m)\le \sum_{m\ge 1}|I_m|\left(\frac{x^\star_{p_m}}{t_m}+\frac{1}{m^2t_m}\right) < \infty,
\end{displaymath}
hence $A \in \mathcal{I}_z$. Similarly, $A-1 \notin \mathcal{I}_z$ since
\begin{displaymath}
\sum_{n \in A-1}z_n
\ge \sum_{m\ge 1}\sum_{n \in I_m \cap 2\mathbf{N}}(z_n^{(m)}-\delta_m)
\gg \sum_{m\ge 1}|I_m|\left(\frac{x^\star_{p_m-1}}{t_m}-\frac{1}{m^2t_m}\right) 
= \infty.
\end{displaymath}
Now, fix $y \in \mathscr{Y}$ such that $\kappa:=\liminf_n y_{n+1}/y_n>0$ and suppose that $\mathcal{I}_z=\mathcal{I}_y$. Then we would have that $A \in \mathcal{I}_y$ and $A-1 \notin \mathcal{I}_y$, which is impossible reasoning as in \eqref{eq:impossible}.

\section{Proof of Theorem \ref{thm:examples}}

%

Let $x^\star$ be the sequence defined in the proof of Theorem \ref{thm:existencefubini}. 
For each $r \in (0,1]$, let $x^{(r)}$ be the sequence defined by $x^{(r)}_n=(x^\star_n)^r$ for all $n \in \mathbf{N}$. Replacing the set of primes $\mathbf{P}$ with $\{\lfloor p_n^{1/r}\rfloor: n \in \mathbf{N}\}$ and reasoning as in the proof of Theorem \ref{thm:existencefubini}, we obtain that $x^{(r)} \in \mathscr{A}$.

To complete the proof, fix reals $r,s$ such that $0<r<s\le 1$. 
Then, it is sufficient to show that the ideals generated by $x^{(r)}$ and $x^{(s)}$ are not isomorphic. To this aim, let $f: \mathbf{N} \to \mathbf{N}$ be a bijection and assume for the sake contradiction that
\begin{equation}\label{eq:contradictionideals}
\forall A\subseteq \mathbf{N}, \quad 
\sum_{n \in A}x^{(r)}_{f(n)}<\infty\,\, \text{ if and only if }\,\,\sum_{n \in A}x^{(s)}_{n}<\infty. 
\end{equation}
Fix $t \in (1,\nicefrac{s}{r})$ and define $S:=\{n \in \mathbf{N}: f(n)>n^t\}$ and $T:=\mathbf{N}\setminus S$. We have $\sum_{n \in S}\frac{1}{f(n)}\le \sum_{n \in S}\frac{1}{n^t}<\infty$. Considering that $f$ is a bijection and the harmonic series is divergent, we obtain that $\sum_{n \in T}\frac{1}{f(n)}=\infty$ (in particular, $T$ is infinite). In addition, since $r<1$ and 
$$
\frac{1}{f(n)} \ll \frac{1}{(f(n) \log(f(n)+1))^r} \le x^{(r)}_{f(n)} \le \frac{1}{f^r(n)},
$$
we get that $\sum_{n \in T}x^{(r)}_{f(n)}=\infty$ and, thanks to \eqref{eq:contradictionideals}, also that $\sum_{n \in T}x^{(s)}_{n}=\infty$. Note that, if $n \in T$, then $f(n) \le n^t$, which implies that
$$
\forall n \in T, \quad \frac{x^{(s)}_n}{x^{(r)}_{f(n)}} \le \frac{1/n^s}{1/(f(n) \log(f(n)+1))^r} \ll \frac{(\log n)^r}{n^{s-tr}},
$$
which has limit $0$ if $n\to \infty$ (and belongs to $T$). In particular, for each $k \in \mathbf{N}$, there exists $n_k \in \mathbf{N}$ such that $x^{(s)}_n/x^{(r)}_{f(n)}\le \nicefrac{1}{k^2}$ for all $n\ge n_k$. Let $(A_k)$ be a sequence of finite subsets of $T$ defined recursively as follows: for each $k \in \mathbf{N}$, let $A_k$ be a finite subset of $T$ such that $\min A_k\ge n_k+ \max A_{k-1}$ and $\sum_{n \in A_k}x^{(r)}_{f(n)} \in (\frac{1}{2},1)$ where, by convention, we assume $\max A_0:=0$ (note that it is really possible to define such sequence). Finally, define $A:=\bigcup_k A_k$ so that we obtain
$$
\sum_{n \in A}x^{(r)}_{f(n)}=\sum_{k\ge 1}\sum_{n \in A_k}x^{(r)}_{f(n)}=\infty 
\,\,\,\text{ and }\,\,\,
\sum_{n \in A}x^{(s)}_{n} \le \sum_{k\ge 1}\sum_{n \in A_k}\frac{x^{(r)}_{f(n)}}{k^2}\le \sum_{k\ge 1}\frac{1}{k^2} <\infty.
$$
This contradicts \eqref{eq:contradictionideals}, concluding the proof.

\section{Concluding Remarks}

We remark that the ideal $\mathcal{I}_x$ defined in the proof above is just (an isomorphic copy of) the Fubini sum $\mathcal{I}_{s}\oplus \mathcal{I}_{t}$, where $s,t \in \mathscr{Z}$ are sequences defined by $s_n=\frac{1}{n}$ and $t_n=\frac{1}{n\log(n+1)}$ for all $n \in \mathbf{N}$. Here, we recall that the Fubini sum of two ideals $\mathcal{I}$ and $\mathcal{J}$ on $\mathbf{N}$ is the ideal $\mathcal{I}\oplus \mathcal{J}$ on $\{0,1\}\times \mathbf{N}$ of all sets $A$ such that $\{n \in \mathbf{N}: (0,n) \in A\} \in \mathcal{I}$ and $\{n \in \mathbf{N}: (1,n) \in A\} \in \mathcal{J}$, cf. e.g. \cite[p. 8]{MR1711328}

Also, some comments are in order about the simplifications. The assumption that the sequence $x$ has positive elements (instead of nonnegative elements) is rather innocuous. Indeed, in the opposite, if $x_n=0$ for infinitely many $n$, then the summable ideal $\mathcal{I}_x$ would be (isomorphic to) the Fubini sum $\mathcal{P}(\mathbf{N}) \oplus \mathcal{I}_y$, for some $y \in \mathscr{X}$. Lastly, also the hypothesis $x \in \mathscr{Z}$ in Question \ref{question:filipczak} (instead of $x \in \mathscr{X}$) has a similar justification. Indeed, if $x \in \mathscr{X} \setminus \mathscr{Z}$, then $\mathcal{I}_x$ would be the Fubini sum $\mathrm{Fin} \oplus \mathcal{I}_y$, for some $y \in \mathscr{Z}$. 

We conclude with two open questions:
\begin{question}
Is it true that for each $x \in \mathscr{Z}$ there are (possibly infinite) sequences $y^1,y^2,\ldots \in \mathscr{Y}$ such that $\mathcal{I}_x=\oplus_i \mathcal{I}_{y^i}$?
\end{question}

\begin{question}
Is it true that $\mathscr{A}$ is an analytic subset of $\mathscr{Z}$?
\end{question}

\bibliographystyle{amsplain}
\bibliography{serie}

%
%
%

%

\end{document}